\documentclass[11pt]{amsart}
\usepackage{amsopn, amsfonts,amsmath}
\usepackage{amssymb, amscd}
\newcommand{\C}{\mbox{$\mathbb C$}}
\newcommand{\q}{\mbox{$\mathbb Q$}}
\newcommand{\Z}{\mbox{$\mathbb Z$}}

\newcommand{\n}{\mbox{$\mathfrak{n}$}}
\newcommand{\ab}{\mbox{$\mathfrak{a}$}}
\newcommand{\m}{\mbox{$\mathfrak{m}$}}
\newcommand{\z}{\mbox{$\mathfrak{Z}$}}

\newcommand{\D}{\mbox{deg$_\mathbb Q$}}
\newcommand{\DA}{\mbox{deg$_{\mathbb Q(\alpha)}$}}

\newcommand{\DAB}{\mbox{deg$_{\mathbb Q(\alpha \beta)}$}}

\newcommand{\dg}{\operatorname{deg}}

\theoremstyle{plain}
\newtheorem{theorem}{Theorem}
\newtheorem{proposition}{Proposition}
\newtheorem{corollary}{Corollary}
\newtheorem{lemma}{Lemma}

\theoremstyle{definition}
\newtheorem{definition}{Definition}
\newtheorem{remark}{Remark}

\begin{document}

\title[Anosov Lie algebras and algebraic units]{Anosov Lie algebras and algebraic units in number fields}

\author{Meera G. Mainkar}

\address{Department of Mathematics, Dartmouth College, Hanover, NH 03755, USA } \email{meera.g.mainkar@dartmouth.edu}


\maketitle
\begin{abstract}
We study nilmanifolds admitting Anosov automorphisms by
applying elementary properties of algebraic units in number
fields to the associated Anosov Lie algebras.
 We identify obstructions to the existence of Anosov Lie
 algebras. The case of $13$-dimensional Anosov Lie algebras is
 worked out as an illustration of the technique. Also, we
 recapture the following known results: (i) Every $7$-dimensional
 Anosov nilmanifold is toral, and (ii) every $8$-dimensional Anosov Lie algebra with $3$ or $5$-dimensional derived
 algebra contains an abelian factor.

\end{abstract}

\section{Introduction}\label{intro}

A diffeomorphism $f$ of a compact differentiable manifold $M$
is said to be {\em Anosov} if there is a continuous invariant
splitting of the tangent bundle $TM = E^+ \oplus E^-$ such that
$df$ expands $E^+$ and contracts $E^-$ exponentially. Such
diffeomorphisms are important in the study of hyperbolic
dynamics. The only known examples of Anosov diffeomorphisms on
compact manifolds are defined on nilmanifolds, or more
generally, infranilmanifolds. We recall that a {\em
nilmanifold} $N/\Gamma$ is a compact quotient of a simply
connected nilpotent Lie group $N$  by a discrete subgroup
$\Gamma \subset N$; an {\em infranilmanifold} is a manifold
finitely covered by a nilmanifold.
 A. Manning (see \cite{Man}) and J. Franks (see \cite{F})  proved that
   any Anosov diffeomorphism on a nilmanifold  $N/\Gamma$ is topologically conjugate to an {\em Anosov
   automorphism}, i.e. a diffeomorphism of $N/\Gamma$ induced
   by a hyperbolic Lie group automorphism of $N$ mapping
   $\Gamma$ to itself. (``Hyperbolic'' means that no
   eigenvalue of the differential of the automorphism is of
   absolute value $1$.) In \cite{S}, S. Smale
raised the problem of classifying the nilmanifolds admitting
Anosov diffeomorphisms, which in view of Manning and Franks's
results reduces to the classification of nilmanifolds with
Anosov automorphisms.

At the level of Lie Algebras, this problem corresponds to the
classification of {\em Anosov} Lie algebras, i.e., the Lie
algebras of nilpotent Lie groups which have nilmanifold
quotients admitting Anosov automorphisms. For previous work in
this direction, see
\cite{S,A-S,Dn,D-M,D,De-Des,La,lw1,lw2,M,MW,P}. In this paper we
consider a number-theoretic approach, which gives some
information regarding the structure of Anosov Lie algebras. We
work out the case of 13 dimensional Lie algebras as an
illustration of this method, and give simpler proofs of some
results obtained in \cite{lw1,lw2}.

To state our main result, we recall a few definitions
introduced in \cite{lw1}.

\begin{definition} \label{abelian}
{\rm  Let $\n$ be a  Lie algebra. An {\em abelian factor} of
$\n$ is an abelian (Lie) ideal $\ab$ of $\n$ such that $\n = \m
\oplus \ab$ for some ideal $\m$ of $\n$. }
\end{definition}

\begin{definition} \label{type}
{\rm Let $\n$ be an $r$-step nilpotent Lie algebra, i.e., the
lower central series $\{C^i(\n)\}$ (defined by $C^0(\n)=\n$ and
$C^i(\n)=[\n,C^{i-1}(\n)]$ for $i \geq 1$) satisfies
$C^{r-1}(\n)\not=0$ and $C^r(\n)=0$. Then the {\em type}
 of $\n$ is the  $r$-tuple of positive
integers $(n_1,\cdots,n_r)$, where $n_i=  \dim C^{i-1}(\n) /
C^{i}(\n)$.}
\end{definition}

We characterize the complexifications of 13-dimensional Anosov
Lie algebras up to isomorphism:

\begin{theorem} \label{thm-main}
Every 13-dimensional real Anosov Lie algebra without an abelian
factor is of type $(9, 4)$. For each such Lie algebra there
exist complex numbers $a, b, c, d$ such that its
complexification has $\mathbb{C}$-linear basis
\[ X_1,X_2,X_3,Y_1,\cdots,Y_6,Z_1,Z_2,W_1,W_2,\]
where the only non-zero brackets of basis elements are
\begin{equation}\label{9,4}
 \begin{array}{ll}
 [X_1, Y_1] = Z_1, & [X_2, Y_2] = W_1  \\

      [X_1, Y_4] = Z_2, & [X_2, Y_5] = W_2 \\

 [X_3, Y_3] = aZ_1 + b W_1, & [X_3, Y_6] = c Z_2 + d W_2.

 \end{array}
 \end{equation}

\end{theorem}

An actual example of a 13-dimensional Anosov Lie algebra was
given in \cite{MW}. It will be seen later that the hypothesis
in Theorem~\ref{thm-main} regarding absence of an abelian
factor is natural. Our method gives  new proofs of the
following results from \cite{lw1,lw2}:
\begin{theorem}\begin{enumerate}
\item Every  7-dimensional Anosov Lie algebra is abelian.
\item Every $8$-dimensional Anosov Lie algebra of type $(5,3)$
    has an abelian factor.
\item There are no $8$-dimensional Anosov Lie algebras of type
    $(3, 3, 2)$.
\end{enumerate}
\end{theorem}

In the last section, we give a necessary condition for a
 Lie algebra of type $(n, 2)$ to be Anosov for $n$ odd.

\vspace{.5cm}

 \noindent {\it Acknowledgements:}
  I am grateful to Prof. J. Lauret and Prof. C. E. Will for their
help.

\vspace{.5cm}

\section{Some Background}

Following \cite{La}, a rational Lie algebra $\n$ of dimension $n$
is said to be {\em Anosov} if it admits a hyperbolic automorphism $\tau$
(i.e. all eigenvalues of $\tau$ have absolute value different
from $1$), and there is a basis of $\n$ with respect to which
the matrix of $\tau$ lies inside $GL(n, \Z)$. We say that a real
Lie algebra is Anosov  if it admits a rational form which is
Anosov. The map $\tau$ will be called  an Anosov automorphism
of the Lie algebra $\n$.

If a nilmanifold $N / \Gamma$ admits an Anosov automorphism,
then it can be seen  that the rational Lie algebra determined
by the lattice $\Gamma$ is Anosov and hence the Lie algebra of
$N$ is Anosov. Therefore, in order to study the nilmanifolds
admitting Anosov automorphisms, we can equivalently study
Anosov Lie algebras.

We recall  \cite[Theorem 3.1]{lw1}: Let $\n$ be a rational Lie
algebra and let $\n = \m \oplus  \ab$ be a Lie direct sum,
where $\ab$ is a  maximal abelian factor
  of $\n$. Then $\n$ is Anosov iff
$\m$ is Anosov and dim $\ab \geq 2$. In view of this, we are
interested in studying  Anosov Lie algebras without an abelian
factor.

Let $\tau$ be an Anosov automorphism of a Lie algebra $\n$.
Then the eigenvalues of $\tau$ are algebraic units. This
follows since by definition there is a distinguished basis with
respect to which the matrix of $\tau$ lies in $GL(n, \Z)$, so
that the characteristic polynomial of $\tau$ is monic with
integer coefficients and $\pm 1$ constant term.
 It follows that  {\em an eigenvalue of an Anosov automorphism of a Lie
algebra is an  algebraic unit with absolute value different
from 1.} This condition will allow us to prove the
non-existence of Anosov Lie algebras (and therefore
nilmanifolds admitting Anosov automorphisms) in some cases.

We summarize the present state of knowledge regarding Anosov
Lie algebras. There are no non-abelian Anosov Lie algebras of dimension
less than 6 (see \cite{W}). There exists an {\em
indecomposable} (not a Lie direct sum of smaller dimensional
Lie algebras) $2$-step Anosov Lie algebra  of
 dimension $n$, for every integer $n \geq 6 $, $n \neq 7$ (see \cite{D-M,MW}.)
For $n  \geq 17$ and $n = 6, 8, 10, 11, 14, 15$, the examples
of indecomposable $2$-step Anosov Lie algebras are given in
\cite{D-M} which are associated to the graphs. For $n = 9, 12$,
see \cite{La}. For $n = 16$ one can modify the constructions of
\cite{D-M} (see \cite{MW}). These Lie algebras are associated
with certain graphs as in \cite{D-M} or obtained by modifying
such examples.  However the existence for $n = 13$ has been
proved by using  properties of the algebraic units (see
\cite{MW}).

\vspace{0.5cm}

\section{Algebraic Units not on the Unit Circle} \label{basic}

Denote by $[E : F]$ the degree of a field extension $F
\subseteq E$. For $\lambda \in E$  let $F(\lambda)$ denote the
smallest subfield of $E$ containing $\lambda$ and $F$. The
degree of an algebraic element $\lambda \in E$  over $F$ is
denoted by $\dg_F(\lambda)$. We recall that $\dg_F(\lambda) =
[F(\lambda) : F]$. We say that $\lambda' \in E$ is a {\em
conjugate} of $\lambda \in E$ over $F$ if $\lambda$ and
$\lambda'$ satisfy the  same minimal polynomial over $F$, and
in this case $\lambda$ and $\lambda'$ are {\it conjugates over
$F$}. We also note that the conjugate elements over $F$ are
conjugate under the action of the  appropriate Galois group.

\begin{remark}\label{rem1}
{\rm We note that
 $\q(\alpha)(\beta) = \q (\alpha, \beta) =
 \q (\alpha \beta)(\alpha) = \q (\alpha \beta)(\beta)$ for all
 non-zero $\alpha, \beta \in \C$}.
\end{remark}

We call $\lambda \in \C$ an {\it algebraic unit } if it
satisfies a monic polynomial with integer coefficients and with
constant term $\pm 1$. Note that $\lambda$ is an algebraic unit
iff both $\lambda$ and $\lambda^{-1}$ are algebraic integers.

  {\bf Throughout this section we will assume that  $\alpha$ and $\beta$ are
algebraic units} and we will derive some properties of
algebraic units which will be used to prove our main results.

 We have
\begin{equation}\label{degree}
\begin{array}{ll}
[\q(\alpha, \beta) : \q(\alpha)] [\q(\alpha) : \q] &= [\q(\alpha,
\beta)  : \q(\alpha \beta)] [\q(\alpha \beta) : \q]\\
 &= [\q(\alpha,
\beta) : \q(\beta)][\q(\beta) : \q].
\end{array}
\end{equation}

\newpage

\begin{lemma}\label{lemma0}
If $\D(\beta) = \D(\alpha
\beta) = \DA(\beta)$ then $\alpha^{\D(\beta)} = \pm  1$.
\end{lemma}

\proof We fix an algebraic closure $\overline{\mathbb{Q}}$ of
$\q$ containing $\alpha$ and $\beta$. Let $\D(\beta) = n$.
 Let $\{\beta = \beta_1, \cdots, \beta_n \}$
  denote   the set of  all conjugates
 of $\beta$ over $\q$ in $\overline{\mathbb{Q}}$ and let $\{\gamma_1 = \alpha \beta, \cdots,
\gamma_n\}$ denote the set of all conjugates of $\alpha \beta$
over $\q$. Note that the set of  all conjugates of $\alpha
\beta$ over $\q$ is same as the set of all conjugates of
$\alpha \beta$ over $\q(\alpha)$  since  $\D(\alpha \beta) =
\DA(\beta)$; see Remark \ref{rem1}. As noted above, the
conjugates of $\alpha\beta$  (or of $\beta$) in $\overline{\q}$
over $\q(\alpha)$ are the conjugates under the action of the
Galois group of $\overline{\q}$ over $\q(\alpha)$. Hence
$\beta, \alpha^{-1}\gamma_2, \ldots, \alpha^{-1} \gamma_n$ are
all the (distinct) conjugates of $\beta$ over $\q(\alpha)$ and
hence over $\q$  since $\D(\beta) = \DA(\beta)$.
  As $\beta$ is an algebraic unit,
 the product $\prod_{i=1}^{n} \alpha^{-1} \gamma_i = \pm 1$.
  But $\prod_{i=1}^{n} \gamma_i = \pm 1$ because $\alpha\beta$ is an
 algebraic unit. Hence $\alpha^{n} = \pm  1$. \hfill$\square$
 \vspace{.2cm}

\begin{corollary}\label{lemma3}
Suppose that $\D(\alpha) = \D(\beta) = \DAB(\alpha)= m$. Then
$$(\alpha \beta)^{m} = \pm 1.$$
\end{corollary}

\proof The proof is immediate if  we set  $\tilde{\alpha} =
\alpha \beta$ and $\tilde{\beta} = \alpha^{-1}$ and use Lemma
\ref{lemma0}.  \hfill$\square$

\vspace{.2cm}

\begin{remark}\label{lemma4}
{\rm From Equation \ref{degree} and Corollary \ref{lemma3}, we
deduce that if  $\D(\alpha) = \D(\beta)$ is coprime to
$\D(\alpha \beta)$,  then $(\alpha \beta)^{\D(\beta)} = \pm
1$.}
\end{remark}

\vspace{.2cm}

\begin{corollary}\label{lemma1}
If $\D(\alpha)$ is coprime to
 $\D(\beta)$ and $\D(\alpha \beta)$ is prime,  then  $$\alpha^{\D(\alpha
 \beta)} = \pm 1 \mbox{ or } \beta^{\D(\alpha
 \beta)} = \pm 1$$
\end{corollary}

\proof If $m = [\q(\alpha) : \q]$ is coprime to $n =
[\q(\beta) : \q]$, then $[\q(\alpha, \beta) : \q(\alpha)] = n$
and $[\q(\alpha, \beta) : \q(\beta)] = m$ (see Equation
\ref{degree}). If $p = [\q(\alpha \beta) : \q]$ is prime, either $p$ does not divide $m$
or $p$ does not divide $n$.  Suppose that  $p$ does not divide $m$, then $p$ divides $n$.
 But $$n  =
[\q(\alpha)(\alpha \beta) : \q(\alpha)] \leq [\q(\alpha \beta) :
\q] = p,$$ and hence $p = n$. This implies that $\D(\beta) =
\D(\alpha \beta) = \DA(\beta)$. From Lemma \ref{lemma0}, we
deduce that $\alpha^{p} = \pm 1$.
\hfill$\square$

\vspace{.2cm}

\begin{corollary}\label{lemma2}
If  $\D(\alpha)$ is coprime to $\D(\beta)$ and if $\D(\alpha
\beta) = \D(\beta)$,  then $$\alpha^{\D(\beta)} = \pm 1.$$
\end{corollary}

\proof Suppose $\D(\alpha \beta)= \D(\beta)$.
Then $$[\q(\beta) : \q] = [\q(\alpha \beta) : \q] = [\q(\alpha,
\beta) : \q(\alpha)]$$  since
 $[\q(\alpha) : \q]$ and $[\q(\beta) : \q]$ are
relatively prime, and $[\q(\beta) : \q]$ is greater than or
equal to $[\q(\alpha, \beta) : \q(\alpha)]$; see Remark
\ref{rem1} and Equation (\ref{degree}). Hence $$\D(\beta) =
\D(\alpha \beta) = \DA(\beta).$$  From Lemma \ref{lemma0}, we
deduce that $\alpha^{\D(\beta)} = \pm 1$. \hfill$\square$

\vspace{.2cm}

\begin{lemma}\label{lemma5}
If $\D(\alpha) < \D(\beta)$ then $\D(\alpha \beta)$
 cannot be  coprime to $\D(\beta)$.

\end{lemma}

\proof Let $\D(\alpha) = m, \D(\beta) = n$ with $m < n$, and let
 $\D(\alpha \beta) =
l$. Suppose that $l$ is coprime to $n$. Then by Equation
(\ref{degree}), $n = [\q(\alpha, \beta) : \q(\alpha \beta)]$
and $l = [\q(\alpha, \beta) : \q(\beta)]$ since
$[\q(\beta)(\alpha \beta) : \q(\beta)] \leq [\q(\alpha \beta):
\q]$. But then $$n = [\q(\alpha \beta)(\alpha) : \q(\alpha
\beta)] \leq m = [\q(\alpha):\q]$$ which is a contradiction.
\hfill$\square$

\vspace{.2cm}

\begin{lemma} \label{(4+6, 3)}

If $\D(\alpha) = 4, \D(\beta) = 6,$ and $\D(\alpha \beta) = 3$,
then $(\alpha \beta)^3 = \pm 1$.
\end{lemma}

\proof  Let $K = \q(\alpha, \beta)$. Then we have
$[K:\q(\alpha)] = 3$ and hence $[K : \q] = 12$. Also
$[\q(\beta) : \q] = 6$ implies that $[K : \q(\beta)] = 2$. Let
$\alpha_1 = \alpha, \alpha_2$ be the conjugates of $\alpha$
over $\q(\beta)$ and let $\beta_1 = \beta, \beta_2$ and
$\beta_3$ denote the conjugates of $\beta$ over $\q(\alpha)$.
Then the conjugates of $\alpha\beta$ over $\q(\alpha)$ (and
hence over $\q$) are $\alpha_1 \beta_1, \alpha_1 \beta_2,
\alpha_1 \beta_3$. Applying the non-trivial Galois automorphism
$\sigma$ of $K$ over $\q(\beta)$ to $\alpha_1 \beta_1$, we
obtain $\alpha_2 \beta_1$. Without loss of generality, we
assume that
\begin{equation} \label{*}
 \alpha_2 \beta_1 =\alpha_1 \beta_2
 \end{equation}
This shows that $\beta_2$ and $\beta_3$ are in $K$ and hence
$K$ over $\q(\alpha)$ is Galois. We pick a Galois automorphism
$\rho$ of $K$ over $\q(\alpha)$ such that $\rho(\beta_1) =
\beta_2, \rho(\beta_2) = \beta_3$ and $\rho(\beta_3) =
\beta_1$.

We claim that $\alpha_2 \in \q(\alpha_1)$. If $\alpha_2$ is not
in $\q(\alpha_1)$, then $K = \q(\alpha_1,  \alpha_2)$. Applying
$\rho$ to $\alpha_2$, we can conclude that  $K$ is a splitting
field of the minimal polynomial of $\alpha$ over $\q$. Its
Galois group, a subgroup of order 12 of the symmetric group
$S_4$, must be the alternating group $A_4$. Note that $K$
contains a splitting field of the minimal polynomial of $\alpha
\beta$ over $\q$, say $M$. In view of the effect of $\sigma$
and $\rho$, the Galois group of $M$ over $\q$ is $S_3$. But
$A_4$ does not surject onto $S_3$ which a contradiction. Hence
$\alpha_2$ is in $\q(\alpha_1)$.

Since $\rho$ is a Galois automorphism of $K$ over
$\q(\alpha_1)$,
 $\rho(\alpha_2) = \alpha_2$ and if we apply $\rho$ to Equation
 \ref{*}, we get $\alpha_2 \beta_2 = \alpha_1 \beta_3$ and
 $\alpha_2 \beta_3 = \alpha_1 \beta_1$.  This gives
 $\frac{\beta_1}{\beta_2} = \frac{\beta_3}{\beta_1}$ and hence
 $$(\alpha_1 \beta_1)^2 = (\alpha_1 \beta_2) (\alpha_1
 \beta_3).$$ Then $(\alpha_1 \beta_1)^3 = (\alpha_1
 \beta_1)(\alpha_1 \beta_2)(\alpha_1 \beta_3) = \pm 1$. \hfill$\square$

\vspace{.2cm}

\begin{lemma}\label{(6+6, 3)}
Let $\D(\alpha) = 6$ and $\beta$ denote a conjugate of $\alpha$
over $\q$, such that $\DAB(\alpha) = 3$, then $(\alpha \beta)^3
= \pm 1$.
\end{lemma}

\proof We claim that $\alpha$  and $\beta$ are not conjugates over
$\q(\alpha \beta)$. Because if $\alpha$ and $\beta$ are conjugates
over $\q(\alpha \beta)$, then the extension $\q(\alpha, \beta)$ over
$\q(\alpha \beta)$ is cyclic Galois extension of degree 3 and there
exists an automorphism  of $\q(\alpha, \beta)$ of order $3$ mapping
$\alpha$ to $\beta$ and fixing $\alpha \beta$, which is impossible.
Let $\{\alpha_1 = \alpha, \alpha_2, \alpha_3\}$ and $\{\beta =
\beta_1, \beta_2, \beta_3\}$ denote the sets of conjugates of
$\alpha$ and $\beta$ over $\q(\alpha \beta)$ respectively. Then it
can be seen that $\alpha, \beta_2^{-1} \alpha \beta$ and
$\beta_3^{-1}\alpha\beta$ are conjugates over $\q(\alpha \beta)$.
Hence $(\alpha_2 \beta_2)(\alpha_3 \beta_3) = (\alpha\beta)^3$.
 We conclude that $$(\alpha \beta)^3 = \prod_{i=1}^{3} \alpha_i
 \prod_{i=1}^{3}\beta_i =\pm 1$$ since $\alpha$ and $\beta$ are
 algebraic units and conjugates over $\q$. \hfill$\square$

 \vspace{.5cm}

\section{$13$-dimensional Anosov Lie algebras: Proof of Theorem~\ref{thm-main}}
\label{13dim}

Although Theorem~\ref{thm-main} holds for Lie algebras with
arbitrary number of steps, for clarity of exposition we write
the proof for 2-step nilpotent Lie algebras. The idea of the
proof generalizes easily to Lie algebras with three or more
steps, although the details are cumbersome.

 We note that, to prove the  existence  of a
$13$-dimensional indecomposable Anosov Lie algebra, our  known
methods (from \cite{D-M} or \cite{La}) are not sufficient.
  We use  the properties of the special algebraic units arising
  from  Anosov automorphisms.  If we consider all possible cases,
  by using Lemmas from Section \ref{basic}, we rule out all cases
  for $13$-dimensional Lie algebra without an abelian factor to be
   Anosov except for the type $(9, 4)$  (see Definition \ref{type}).

We recall   \cite[Proposition~2.1]{lw1} and state it for
$2$-step Anosov Lie algebras which will be used throughout this
paper.

\begin{proposition}\label{LW}
 Let $\n$ be a $2$-step Anosov Lie algebra. Then there exists
  a hyperbolic semisimple automorphism $\tau$ of $\n$ and a
   vector space decomposition $\n = \n_1 \oplus [\n, \n]$ such that
   $\tau(\n_1) = \n_1$ and the characteristic polynomial $f$ (resp. $g$) of
    the restriction of $\tau$ to $\n_1$ (resp. to $[\n, \n]$) is with integer
    coefficients with constant term $\pm 1$.
\end{proposition}

 Let $\n$ be a $2$-step Anosov Lie algebra and let
 $\tau$, $\n_1$, $f$ and $g$ be chosen as in  Proposition \ref{LW}.
    We note that  none of the roots of $f$
and $g$ are of modulus equal to $1$. All the roots of $g$ are
of the form $\alpha \beta$ where $\alpha$ and $\beta$ are roots
of $f$. Moreover, $f$ and $g$ are monic polynomials with
integer coefficients with  constant term $\pm 1$. Hence the
roots of $f$ and $g$ are algebraic units. We may assume that
the constant term of $f$ and $g$ is $1$ by considering the
square of an automorphism if required.

\begin{lemma}\label{abelianfactor}
Let $\n$ be a $2$-step Anosov Lie algebra with $\n = \n_1
\oplus [\n, \n]$ and $\tau$ be a semisimple hyperbolic
automorphism of $\n$ such that $\tau(\n_1) = \n_1$. Let  $f$
(resp. $g$) denote the characteristic polynomial  of
$\tau|_{\n_1}$ (resp. $\tau|_{[\n, \n]}$). If there exists a
root $\alpha$ of $f$ such that for all roots $\beta$ of $f$ one
has $g(\alpha \beta) \neq 0$ then $\n$ has an abelian factor.
\end{lemma}

\proof Since $\tau$ is a semisimple automorphism, there exists
a basis of $\n_1^{\mathbb{C}}$ (complexification of $\n_1$)
consisting of the eigenvectors corresponding to the eigenvalues
of the restriction of $\tau$ on $\n_1$.
 Let $X$ be a  eigenvector corresponding to the eigenvalue
 $\alpha$. Then $[X, Y] = 0$ for all $Y \in \n_1^{\mathbb{C}}$ by
 our hypothesis. Hence $X$ belongs to the center of
 $\n^{\mathbb{C}}$, $\z(\n^{\mathbb{C}})$. In particular $\z(\n)
 \cap [\n, \n] \neq \z(\n)$, since $\z(\n^{\mathbb{C}}) \cap
 [\n^{\mathbb{C}}, \n^{\mathbb{C}}] = (\z(\n) \cap [\n,
 \n])^{\mathbb{C}}$. This implies that $\n$ has an abelian factor.
 \hfill$\square$

\vspace{1cm}

 Let $\n$ be a $2$-step $13$-dimensional Anosov Lie algebra
without an abelian factor of type $(n_1,
 n_2)$. Then $n_1 \geq 5$ since $\n$ is 13-dimensional and $n_2 \geq 2$
  (equivalently $n_1 \leq 11$) because an Anosov automorphism
   cannot act as 1 or -1 on a one-dimensional center. We
 consider all the cases and we will show that $\n$ must be
 of the type $(9, 4)$. We choose $\tau$, $\n_1$, $f$ and $g$ as
 in Proposition \ref{LW}.
 \vspace{0.5cm}

 \noindent {\bf Case $(11, 2)$}.
 We note that $g$ is irreducible over $\Z$ since the degree of $g$ is $2$
and $g$ does not have an eigenvalue of modulus $1$.
 If $f$ is irreducible over $\Z$ then there exist algebraic units,
 $\alpha$ and $\beta$, such that $\D(\alpha) = \D(\beta) = 11$ and
$\D(\alpha \beta) = 2$. This is not possible by Remark
\ref{lemma4}.
 If $f$ is reducible over $\Z$, then there exists an odd degree
irreducible factor
 of $f$ over $\Z$, say $h$.  The degree of $h$ could be 3, 5, 7 or 9.
  We will prove that the degree of $h$ must be 3.

   Suppose that the degree of $h$ is 5. If $f = h h'$  where $h'$ is an irreducible polynomial of degree 6 over $\Z$, then by using
   Remark \ref{lemma4} and Corollary \ref{lemma1}, we can conclude that  $g(\alpha \beta) \neq 0$ for all roots $\alpha$ of $h$ and
    for all roots $\beta$ of $f$. This is a contradiction (see Lemma \ref{abelianfactor}) because we are assuming that $\n$ is without an abelian factor. If $h'$ is reducible over $\Z$ and is a product of two irreducible polynomials each of degree 3 over $\Z$, then also we conclude that
     $g(\alpha \beta) \neq 0$ for all roots $\alpha$ of $h$ and
    for all roots $\beta$ of $f$ by  Remark \ref{lemma4} and Corollary \ref{lemma1}.  By considering all possibilities for  the degrees of 
     the irreducible factors of $h'$, and using Remark \ref{lemma4}, Corollary \ref{lemma1} and Lemma \ref{abelianfactor}, we can prove that the degree of $h$ cannot be equal to 5.

  Similarly we can see that  the degree of $h$ must be
$3$ and that there exists an irreducible factor of $f$ over
$\Z$, say $h'$, of degree $6$; see Corollary \ref{lemma1},
Remark \ref{lemma4} and Lemma \ref{abelianfactor}.
 Let  $f =
 h h' h''$ where $h''$ is an  irreducible polynomial of degree 2 over
$\Z$.
 The product $\alpha \beta$, where $\alpha$ is a root of $h$ and
$\beta$ is a root of $h''$, cannot occur as a root of $g$ by Corollary
\ref{lemma1}. If $\alpha$ is a root of $h'$ and $\beta$ is a root of
$h''$, then we have $\D(\alpha) = 6$ and $\D(\beta) = 2$. In this case
 $\D(\alpha \beta) \neq 2$ as $[\q(\alpha , \beta) : \q(\alpha \beta)]
 \leq 2$ and $[\q(\alpha , \beta) : \q] \geq 6$; see Remark \ref{rem1} and Equation (\ref{degree}).
 Thus there exists an abelian factor (see Lemma \ref{abelianfactor}) which is a contradiction.
 Hence there does not exist an Anosov Lie algebra without an abelian
factor of type $(11, 2)$.

\vspace{0.5cm}

 \noindent {\bf Case $(10, 3)$}. Also  in this case $g$
must be irreducible over $\Z$.  Moreover,
    $f$ cannot be  irreducible by Remark \ref{lemma4}.
    If $f$ has two irreducible
factors, say $f_1$ and $f_2$ over $\Z$ of degree 2 and 8
respectively, and if $\alpha$ is a root of $f_2$, then by using
Remark \ref{lemma4}  it can be seen that for all roots $\beta$
of $f_2$, $g(\alpha \beta) \neq 0$. Also, by using Lemma
\ref{lemma5}, we deduce that for all roots $\beta$ of $f_1$,
$g(\alpha \beta) \neq 0$. This is not possible since $\n$ has
no abelian factor (see Lemma \ref{abelianfactor}). Similarly,
 considering all possibilities for the factorization of $f$ and
   using  Remark \ref{lemma4},  Lemma \ref{lemma5}, Corollary \ref{lemma1} and Lemma \ref{abelianfactor},
   it can be seen that either $f = f_1 f_2$ where
$f_1$ and $f_2$ are irreducible factors of $f$  over $\Z$
 of degree $4$ and $6$ respectively; or $f = h_1 h_2 h_3$ where $h_1$,
$h_2$ and $h_3$ are irreducible factors of $f$ over $\Z$ such
that deg $h_1$ = deg $h_2  = 2$ and deg $h_3  = 6$.

Suppose that  $f = f_1 f_2$ where $f_1$ and $f_2$ are
irreducible factors of $f$  over $\Z$ of degree $4$ and $6$
respectively. Since $\n$ does not have an abelian factor, there
exist roots $\alpha$ of $f_1$ and $\beta$ of $f_2$ such that
$g(\alpha \beta) = 0$ by Lemma \ref{abelianfactor} and Remark
\ref{lemma4}. But this is not possible by Lemma \ref{(4+6, 3)}.

For the other case, we first prove the following lemma:

 \begin{lemma}\label{(2+2+6,3)}
 Let $\beta_1, \beta_2, \beta_3$ denote the distinct roots of an
 irreducible polynomial of degree 3 over $\Z$ with $\pm 1$
 constant term.
 Let $\lambda$ and $\mu$ be real algebraic units such that $$\D(\lambda)
=
 2 = \D(\mu).$$ If $\{ \mu
\beta_1, \mu \beta_2, \mu \beta_3, \mu^{-1} \beta_1,
 \mu^{-1} \beta_2, \mu^{-1} \beta_3 \} = \{ \lambda \beta_1,
 \lambda \beta_2, \lambda \beta_3, \lambda^{-1} \beta_1,
 \lambda^{-1} \beta_2, \lambda^{-1} \beta_3 \}$, then either
 $\mu = \lambda$ or $\mu = \lambda^{-1}$.
 \end{lemma}

\proof If $\mu \neq \lambda$ and  $\mu \beta_1 = \lambda
\beta_i$ for some $i \in \{2, 3\}$, then $\beta_1 \beta_i^{-1}
= \lambda \mu^{-1}$. But $\D(\lambda \mu^{-1}) \in \{2, 4\}$
and
 $\D(\beta_1) = \D(\beta_i^{-1}) = 3$.
 By Remark \ref{lemma4}, $\lambda^3 = \mu^3$. This is a
 contradiction because $\mu$ and $\lambda$ are reals and $\mu \neq
 \lambda$.
 Similarly we see that if
 $\mu \neq \lambda^{-1}$ then $\mu \beta_1 \neq \lambda^{-1} \beta_i$
for all $i \in \{2, 3\}$.
 Hence $\mu = \lambda$ or $\mu = \lambda^{-1}$. \hfill$\square$

Now suppose that $f = h_1 h_2 h_3$ where $h_1$, $h_2$ and $h_3$
are irreducible factors of $f$ over $\Z$ such that deg $h_1$ =
deg $h_2  = 2$ and deg $h_3  = 6$.
 Since $h_1$ and $h_2$ are monic polynomials of degree $2$ with unit
constant term, we may assume that $\lambda$ and $\lambda^{-1}$
are roots of $h_1$, $\mu$ and $\mu^{-1}$ are roots of $h_2$. We
note that $\lambda$ and $\mu$ are reals since they are of
degree 2 algebraic units which are not on the unit circle.  Let
$\{\alpha_1, \ldots, \alpha_6\}$ denote the set of roots of
$h_3$ and let
 $\{\beta_1, \beta_2, \beta_3\}$ denote the set of roots of $g$.
 Since $\n$ is without an abelian factor, there exist $\alpha_i$
 and $\alpha_j$ such that $\mu \alpha_i$ and $\lambda \alpha_j$ are roots of $g$, $1 \leq i, j, \leq 6$.
   Hence the set of roots of $h_3$ is  $\{ \mu
\beta_1, \mu \beta_2, \mu \beta_3, \mu^{-1} \beta_1,
 \mu^{-1} \beta_2, \mu^{-1} \beta_3 \}$ which is the same as
 $\{ \lambda \beta_1,
 \lambda \beta_2, \lambda \beta_3, \lambda^{-1} \beta_1,
 \lambda^{-1} \beta_2, \lambda^{-1} \beta_3 \}$. By  Lemma
 \ref{(2+2+6,3)}, we conclude that either $\mu= \lambda$ or
 $\mu = \lambda^{-1}$. Without loss of generality we assume
 that $\mu = \lambda$.

We recall some notations which were introduced at the beginning of
this section. $\n = \n_1 \oplus [\n, \n]$ and $\tau$ is the
semisimple Anosov automorphism of $\n$ such that $\tau(\n_1) = \n_1$
and $f$ is the characteristic polynomial of $\tau|_{\n_1}$.  Now let
$\n^{\mathbb{C}}$ (resp. $\n_1^{\mathbb{C}}$)
  denote the complexification of $\n$ (resp. $\n_1$).
 Let $\{X_1, X_2, Y_1, Y_2,$ $Z_1, Z_2, \cdots, Z_6 \}$ denote a basis
 of $\n_1^{\mathbb{C}}$  such that
 $ X_1, X_2, Y_1, Y_2, Z_1, Z_2, \cdots, Z_6$ are eigenvectors of
 $\tau|_{\n_1}$ corresponding to the eigenvalues $\lambda,
\lambda^{-1},
 \lambda , \lambda^{-1}, \lambda \beta_1, \lambda \beta_2, \lambda
\beta_3,
 \lambda^{-1}  \beta_1, \lambda^{-1}  \beta_2,~~~~~~~~~~~~~\lambda^{-1}
\beta_3$ respectively.
 Then $X_1 - c Y_1$ is in the center of $\n^{\mathbb{C}}$ for some $c
\in \C$, which would generate an abelian factor.
 Hence there does not exist an Anosov Lie algebra without an abelian
factor of type $(10, 3)$.

\vspace{0.5cm}

\noindent {\bf Case $(8, 5)$.} By  Corollary \ref{lemma1},
Corollary \ref{lemma2}, Remark \ref{lemma4} and Lemma
\ref{lemma5}, we can see that $f = f_1 f_2$ and $g = g_1 g_2$
such that $f_1, f_2, g_1$ $g_2$ are irreducible polynomials
over $\Z$ of degree $2, 6, 2$ and $3$ respectively. There exist
roots $\alpha$  and  $\beta$ of $f$ such that $g_1(\alpha
\beta)= 0$. We note that $\alpha$ and $\beta$ cannot both be
roots of $f_1$. If $f_1(\alpha) = 0$ and $f_2(\beta) = 0$, then
$[\q(\alpha\beta)(\alpha) : \q(\alpha \beta)] \in \{1,2\}$,
which gives a contradiction since $[\q(\alpha \beta): \q] = 2$
(see Equation (\ref{degree})). Hence $\alpha$ and $\beta$ both
should be roots of $f_2$. We note that
$[\q(\alpha)(\alpha\beta): \q(\alpha)] \in \{1, 2\}$. If
$[\q(\alpha)(\alpha \beta):\q(\alpha)]=1$, then $[\q(\alpha,
\beta): \q(\alpha \beta)] = 3$. By Lemma  \ref{(6+6, 3)}, we
conclude that $(\alpha \beta)^3 = \pm 1$ which is a
contradiction. If $[\q(\alpha)(\alpha \beta):\q(\alpha)] = 2$,
then $[\q(\alpha, \beta) : \q(\alpha \beta)] = 6$ which is  a
contradiction  by Corollary \ref{lemma3}.

\vspace{0.5cm}

\noindent {\bf Case $(7, 6)$}.  By  Corollary \ref{lemma1},
 Corollary \ref{lemma2}, Remark \ref{lemma4} and  Lemma
\ref{lemma5} we see that either (i) $f$ factors as a product of
irreducible polynomials $f_1, f_2$ and $f_3$ over $\Z$ such that
 deg $f_1 =$ deg $f_2$ = 2 and deg $f_3 = 3$ and $g$ is an irreducible
polynomial of degree $6$ over $\Z$, or (ii) $f = h_1 h_2$ where $h_1$ is
a degree $3$ irreducible polynomial and $h_2$ is
a degree $4$ irreducible polynomial over $\Z$ and $g$ is an
irreducible polynomial over $\Z$.

In (i), let $\lambda$ and $\mu$ denote the roots of $f_1$ and $f_2$
respectively. Let $\alpha_1, \alpha_2, \alpha_3$ denote the roots of
$f_3$. Since $\n$ is without an abelian factor, $\lambda \alpha_i$
and $\mu \alpha_j$ must be roots of $g$. Hence $$\{\lambda \alpha_1,
\lambda \alpha_2,\lambda \alpha_3,\lambda \alpha_1^{-1},\lambda
\alpha_2^{-1},\lambda \alpha_3^{-1}\} = \{\mu \beta_1, \mu
\beta_2,\mu \beta_3,\mu \beta_1^{-1},\mu \beta_2^{-1},\mu
\beta_3^{-1}\}.$$ By Lemma \ref{(2+2+6,3)}, $\mu = \lambda$ or $\mu
= \lambda^{-1}$. Suppose $\mu = \lambda$. We recall that $f$ is the
characteristic polynomial of $\tau|_{\n_1}$ where $\tau$ is a
semisimple Anosov automorphism of $\n$ such that $\tau(\n_1) = \n_1$
and $\n = \n_1 \oplus [\n, \n]$. Let $\{X_1, X_2, Y_1, Y_2, Z_1,
Z_2, Z_3\}$ denote a basis of $\n_1^{\mathbb{C}}$ (complexification
of $\n_1$)  such that $X_1, X_2, Y_1, Y_2, Z_1, Z_2, Z_3$ are the
 eigenvectors of $\tau|_{\n_1}$ corresponding to the eigenvalues
 $\lambda, \lambda^{-1}, \mu, \mu^{-1}, \alpha_1, \alpha_2,
 \alpha_3$ respectively. Then $X_1 - a Y_1$ is in the center of
 $\n^{\mathbb{C}}$ which is a contradiction to our assumption that
 $\n$ has no abelian factor.

 In  (ii), if $\alpha$ and $\beta$ are  roots of $f$ such that
 $g(\alpha \beta) = 0$, then $\alpha$ and $\beta$ cannot be both
 roots of $h_1$ since $\D(\alpha \beta) = 6$. We can assume that
 $h_1(\alpha) = 0 = h_2(\beta)$. Then by applying Lemma \ref{(4+6, 3)} to $\tilde{\alpha} = \beta^{-1}$
  and $\tilde{\beta}= \alpha \beta$, we see that $\alpha^3 = \pm 1$,
  which is a contradiction. Hence we conclude that there is no
  Anosov Lie algebra of type $(7,6)$ without an abelian factor.

\vspace{0.5cm}

\noindent {\bf Case $(6, 7)$.} By Remark \ref{lemma4}, Lemma
\ref{lemma5} and Lemma \ref{abelianfactor}, we can see that the
only possibility in this case is the following: $f$ is
irreducible over $\Z$, and $g = g_1 g_2$ where $g_1$ and $g_2$
are irreducible over $\Z$ of degree $3$
  and $4$ respectively.
In this case we will prove that there does not exist
 a root $\alpha \beta$ of $g_2$ such that $\alpha$ and $\beta$ are
roots of $f$.
 Suppose that $\alpha$ and $\beta$ are roots of $f$ such that $\alpha
\beta$ is a root of $g_2$. Let $n = [\q(\alpha \beta)(\alpha) :
\q(\alpha \beta)]$.
 Then we note that $n \in \{3, 6\}$ (by Remark \ref{rem1} and Equation
 (\ref{degree})). If $n = 6$, we get a contradiction by Lemma \ref{lemma3}.
 Hence $n = 3$. But then by Lemma \ref{(6+6, 3)}, $(\alpha
 \beta)^3 = 1$, which is a contradiction.

\vspace{0.5cm}

\noindent{\bf Case $(5, 8)$.} By  Remark \ref{lemma4} and Lemma
\ref{lemma5} it can be seen that there does not exist an Anosov
Lie algebra without an abelian factor and of type $(5, 8)$.

\vspace{0.5cm}

 \noindent{\bf Case $(9, 4)$.} By  Remark
\ref{lemma4} and Lemma \ref{lemma5}, we see that there are two
possibilities in this case: (i) $f = f_1 f_2$ where $f_1$ and
$f_2$ are irreducible polynomials
 over $\Z$ of degree $3$ and $6$ respectively, and $g$
 is irreducible over $\Z$. (ii) $f = h_1 h_1$ and $g = g_1 g_2$ such
that $h_1, h_2, g_1$ and $g_2$ are irreducible polynomials over $\Z$
 and deg $h_1 = 3,$ deg $h_2 = 6$ and deg $g_1$ = deg $g_2$ = $2$.

In case (i), there exists roots $\alpha$ of $f_1$ and $\beta$
of $f_2$ such that $g(\alpha \beta) = 0$ (by Remark
\ref{lemma4} and Lemma \ref{abelianfactor}). By applying Lemma
\ref{(4+6, 3)} to $\tilde{\alpha} = \alpha \beta$ and
$\tilde{\beta}= \beta^{-1}$, we conclude that $\alpha^3 = \pm
1$, a contradiction. Hence case (i) is not possible.

In case (ii), we will prove that there exists an indecomposable
Anosov Lie algebra i.e. an Anosov Lie algebra which cannot be
written as a direct product of
 two proper ideals (see \cite{D-M}).

 Suppose that $f = h_1 h_1$ and $g = g_1 g_2$ such that $h_1,
h_2, g_1$ and $g_2$ are irreducible polynomials over $\Z$
 and deg $h_1 = 3,$ deg $h_2 = 6$ and deg $g_1$ = deg $g_2$ = $2$.
 Let $\{ \alpha_1, \alpha_2, \alpha_3 \}, \{ \lambda, \lambda^{-1} \}$
and $\{ \mu, \mu^{-1} \}$ denote the sets of roots of $h_1, g_1$ and
$g_2$ respectively.
 Since there is no abelian factor,   there exist $\alpha$ and $\beta$,
roots of $h_1$ and $h_2$ respectively, such that $\alpha \beta
= \lambda$. Similarly there exist $\alpha'$ and $\beta'$, roots
of $h_1$ and $h_2$ respectively, such that $\alpha'\beta' =
\mu$. Then it can be seen that the set of roots of $h_2$ is
given by
\begin{align*}
\{ \lambda \alpha_1^{-1}, \lambda \alpha_2^{-1},
  \lambda \alpha_3^{-1},\lambda^{-1} \alpha_1^{-1},  \lambda^{-1}
  \alpha_2^{-1},
 \lambda^{-1} \alpha_3^{-1} \} = \\
   \{ \mu \alpha_1^{-1}, \mu \alpha_2^{-1},
\mu \alpha_3^{-1},\mu^{-1} \alpha_1^{-1},  \mu^{-1}
  \alpha_2^{-1}, \mu^{-1} \alpha_3^{-1}\}.
  \end{align*}
   Hence by Lemma \ref{(2+2+6,3)},  $\mu = \lambda$ or $\mu
= \lambda^{-1}$. Without loss of
 generality we assume that $\mu = \lambda$.

Let $\n^{\mathbb{C}}$ ($\n_1^{\mathbb{C}}$)
 denote the complexification of $\n$ ($\n_1$ respectively).
  Let $X_1, X_2,$ and $X_3$ denote the eigenvectors (in
$\n_1^{\mathbb{C}}$) of $\tau|_{\n_1}$ corresponding to the eigenvalues
$\alpha_1, \alpha_2,$ and $\alpha_3$ respectively.
 Let $Y_1, Y_2, \ldots, Y_6$ denote the eigenvectors in
$\n_1^{\mathbb{C}}$ corresponding to the eigenvalues $
  \lambda \alpha_1^{-1}, \lambda \alpha_2^{-1},$ $
  \lambda \alpha_3^{-1},\lambda^{-1} \alpha_1^{-1},  $ $\lambda^{-1}
  \alpha_2^{-1},
 \lambda^{-1} \alpha_3^{-1}$ respectively. Let $Z_1, W_1$ be linearly independent eigenvectors
corresponding to an eigenvalue $\lambda$ and let $Z_2, W_2$ be
 linearly independent eigenvectors
corresponding to an eigenvalue $\lambda^{-1}$. We may assume that
the Lie brackets are given by the following relations:
 \begin{equation}
 \begin{array}{ll}
 [X_1, Y_1] = Z_1, & [X_2, Y_2] = W_1  \\

      [X_1, Y_4] = Z_2, & [X_2, Y_5] = W_2 \\

 [X_3, Y_3] = aZ_1 + b W_1, & [X_3, Y_6] = c Z_2 + d W_2

 \end{array}
 \end{equation}

  where
$a, b, c, d \in \C$ and other commutators in the generators
vanish. This completes the proof of the main Theorem
\ref{thm-main} for $2$-step Anosov Lie algebras. As noted in
the beginning of this section, the same idea leads to a proof
of the general theorem for Lie algebras with arbitrary number
of steps. In the general proof, Proposition~\ref{LW} has to be
replaced by  \cite[Proposition~2.1]{lw1}

\begin{remark}
For an actual example of a $13$-dimensional Anosov Lie algebra,
 see \cite{MW}.
\end{remark}

\section{More nonexistence results}\label{7dim}
As noted in the introduction (\S\ref{intro}), it is known that
there does not exist a non-toral 7-dimensional Anosov Lie
algebra (see \cite{lw2}). Here we give an alternative proof of
this fact using the results from  \S\ref{basic}. First we note
that a $7$-dimensional Anosov Lie algebra has to be $2$-step
and of the type $(4, 3)$ or $(5, 2)$ (see \cite[Proposition
2.3]{lw1}).

As in the 13-dimensional case (\S\ref{13dim}), following
Proposition \ref{LW}, we choose  a decomposition of a
7-dimensional Anosov Lie algebra $\n$ as $\n = \n_1 \oplus [\n,
\n]$ and a hyperbolic semisimple automorphism $\tau$ such that
$\tau(\n_1) = \n_1$  such that the characteristic polynomial
$f$ (resp. $g$) of the restriction of $\tau$ to $\n_1$ (resp.
to $[\n, \n]$) is with integer coefficients and constant term $\pm 1$. Then we note that
the roots of $g$ are certain products of the roots of $f$. The
roots of $f$ and $g$ are algebraic units and none of them is of
absolute value 1.

\vspace{0.5cm}

 \noindent {\bf Case $(4, 3)$.} In this case, $g$
is of degree 3 and hence has to be irreducible. Now either $f$
is irreducible or $f$ is a product of two degree 2 irreducible
polynomials. By using Remark \ref{lemma4}, we can see that both
cases are impossible. \vspace{0.5cm}

\noindent {\bf Case $(5, 2)$.}  In this case also $g$ is
irreducible and is of degree 2.  It can be seen (using Remark
\ref{lemma4}) that $f$ is reducible and  $f$ is a product of
two irreducible polynomials $f_1$ and $f_2$ such that $f_1$ is
of degree 2 and $f_2$ is of degree 3. Now if $\alpha$ is a root
of $f_1$ and $\beta$ is a root of $f_2$ then $\alpha \beta$
will not occur as a root of $g$ by  Lemma \ref{lemma5}.
 Moreover, the product of two roots of $f_1$ is $\pm 1$ and the product
 of two roots of $f_2$ cannot be of degree 2 by Remark \ref{lemma4}.
 Hence this case is also impossible.

 In a similar way, we can use Corollary \ref{lemma2} and Remark \ref{lemma4}  to give a simpler proof to
 show that there are no Anosov Lie
 algebras of type  $(5,3)$ with no abelian factor and there are no Anosov Lie algebras of type
 $(3, 3, 2)$ with no abelian factor; this  was first proved in \cite{lw1}.

\section{Anosov Lie algebras of type $(n, 2)$, $n$ odd}\label{n}

 In this section we study $2$-step Anosov  Lie algebra of type
 $(n, 2)$  (See Definition \ref{type}), where $n$ is odd. Let
 $\n$ be a 2-step Anosov Lie algebra of type $(n,2)$.
  Consider a vector space decomposition $\n = \n_1
 \oplus [\n, \n]$  of $\n$. We say that $\n_1$ is {\it decomposable}
 if $\n_1 = V_1 \oplus V_2$ where $V_1$ and $V_2$ are nontrivial subspaces of
 $\n_1$ such that $[V_1, V_2] = \{0\}$.  We will prove that if a $2$-step nilpotent Lie algebra $\n$ of
 type $(n, 2)$, with $n$ odd,  is Anosov then there exists
 $\n_1$  (such  that $\n = \n_1
 \oplus [\n, \n]$) which is decomposable.

From now onwards we will assume that $\n$ is a $2$-step nilpotent
Lie algebra of type $(n, 2)$ such that $n$ is odd.
 Suppose that $\n$ is an Anosov Lie algebra.
 Now since $\n$ is Anosov 2-step, we choose $\tau, \n_1, f$ and $g$ as in
  Proposition \ref{LW}.  We note that the roots of $g$ are real since deg $g$ = 2
  and $\tau$ is hyperbolic.

Let $\alpha$ denote a real  root of $f$ such that $\D(\alpha)$ is
odd. Since $n$ is odd, $\alpha$ exists. Let $\mu$ denote a
 root of $g$.
  Since $\D(\mu) = 2$ and $\mu$ is an algebraic unit, $\mu^{-1}$ is
   also a root of $g$.
  Let $S$ denote the set of
  eigenvalues of $\tau|_{\n_1}$ which are contained
 in $\{\mu^{2n} \alpha, \mu^{2n+1} \alpha^{-1} : n \in \Z\}$.
  Let $V_1$ denote the sum of
 the eigenspaces corresponding to all eigenvalues in $S$.
   We note that $V_1$ is a (real) subspace of $\n_1$ and $\tau(V_1) = V_1$.
    Since $\tau|_{\n_1}$ is semisimple, there exists a (real) subspace  $V_2$ of
     $\n_1$ such that $\tau(V_2) = V_2$ and $\n_1 = V_1 \oplus V_2$.
   Since $\alpha$ is an eigenvalue of  $\tau|_{\n_1}$, $V_1$ is
   nontrivial. Let $\beta$ be a conjugate of $\alpha$ over $\q$, $\beta \neq \alpha$.
    Then $\D(\beta) = \D(\alpha)$ which is odd. We note that
    $\beta$ is an eigenvalue of $\tau|_{\n_1}$ which is not
    contained in $S$ because $\D((\mu)^{l}) = 2$ for all $l \neq
    0$ (see Remark \ref{lemma4}).
    Hence $V_2$ is nontrivial. Since the
   eigenvalues of $\tau|_{\n_2}$ are of the form $\mu$ and
   $\mu^{-1}$, $[V_1 , V_2] = \{0\}$. This proves that $\n_1$ is
   decomposable. Hence we have proved the following:

   \begin{proposition}
  If $\n$ is a 2-step Anosov nilpotent Lie algebra of an odd dimension
  such that dim $[\n, \n]$ = 2, then there exists a decomposable $\n_1$, a vector space complement of $[\n, \n]$
   in $\n$.
 \end{proposition}

\end{document}